\documentclass[11pt,a4paper]{amsart}
\usepackage{amsmath,amsfonts,amsthm,amssymb,paralist}
\usepackage[margin=3cm,twoside]{geometry}
\usepackage[utf8]{inputenc}

\title[A note on CR mappings of positive codimension]{A note on CR mappings of positive codimension}

\author{Jean-Charles Suny\'e}

\address{Universit\'e de Rouen, CNRS, Laboratoire de Math\'ematiques Rapha\"el Salem, Avenue de l'Universit\'e, B.P. 12, 76801 Saint Etienne du Rouvray, France}

\email{jean-charles.sunye@etu.univ-rouen.fr}

\thanks{The author was partially supported by the Amadeus program of the ``Partenariat Hubert Curien''.}

\subjclass[2000]{32V10, 32H02, 32V20} \keywords{CR mapping, Artin approximation theorem.}

\newtheorem{Thm}{Theorem}[section]
\newtheorem{Cor}[Thm]{Corollary}
\newtheorem{Pro}[Thm]{Proposition}
\newtheorem{Lem}[Thm]{Lemma}

\theoremstyle{definition}
\newtheorem{Def}[Thm]{Definition}

\theoremstyle{remark}
\newtheorem{Rem}[Thm]{Remark}
\newtheorem{Exa}[Thm]{Example}

\numberwithin{equation}{section}

\def\bl{\begin{Lem}}
\def\el{\end{Lem}}
\def\bp{\begin{Pro}}
\def\ep{\end{Pro}}
\def\bt{\begin{Thm}}
\def\et{\end{Thm}}
\def\bc{\begin{Cor}}
\def\ec{\end{Cor}}
\def\bd{\begin{Def}}
\def\ed{\end{Def}}
\def\be{\begin{Exa}}
\def\ee{\end{Exa}}
\def\bpf{\begin{proof}}
\def\epf{\end{proof}}
\def\ben{\begin{enumerate}}
\def\een{\end{enumerate}}
\def\beq{\begin{equation}}
\def\eeq{\end{equation}}
\def\beqar{\begin{eqnarray}}
\def\eeqar{\end{eqnarray}}
\def\brm{\begin{Rem}}
\def\erm{\end{Rem}}

\def\ov{\overline}

\def\zb{\bar z}

\def\z{\zeta}
\def\phi{\varphi}

\def\T{\Theta}

\def\a{\alpha}

\def\rp{{\rho^\p}}
\def\C{\mathbb C}

\def\N{\mathbb N}
\def\R{\mathbb R}

\def\d{\partial}

\def\CR{\text{\rm CR }}

\def\Mp{M^\p}
\def\Np{N^\p}

\def\zp{{z^\p}}

\def\wp{{w^\p}}

\def\pb{\bar{p}}
\def\calT{\mathcal{T}}
\def\calW{\mathcal{W}}
\def\Htilde{\widetilde{H}}
\def\Fp{F^\p}
\newcommand{\p}{\prime}

\newcommand{\Cinf}{\mathcal{C}^{\infty}}

\newcommand{\nequiv}{{\equiv \!\!\!\!\!\!  / \,\,}}

\newcommand{\mc}[1]{\mathcal #1}

\newcommand{\deri}[2]{\frac{\d #2}{\d #1}}
\newcommand{\mderi}[3]{\frac{\d^{#2} #3}{\d {#1}^{#2}}}

\begin{document}

\begin{abstract}
We prove the following Artin type approximation theorem for smooth CR mappings: given $M\subset \C^N$ a connected real-analytic CR submanifold that is minimal at some point, $\Mp\subset \C^{\Np}$ a real-analytic subset, and $H\colon M\to \Mp$ a $\Cinf$-smooth CR mapping, there exists a dense open subset $\mc{O}\subset M$ such that for any $q\in \mc{O}$ and any positive integer $k$ there exists a germ at $q$ of a real-analytic CR mapping $H^k\colon (M,q)\to \Mp$ whose $k$-jet at $q$ agrees with that of $H$ up to order $k$.
\end{abstract}

\maketitle

\section{Introduction}\label{s:intro}

Given germs of real-analytic submanifolds $M$ and $\Mp$ embedded in complex spaces, a fundamental question is to decide whether the formal equivalence of $M$ and $\Mp$ implies their biholomorphic equivalence. While this need not be in general the case in view of a well known example due to Moser-Webster \cite{MW83} (see also \cite{Hua04, HY09}), recent results due to Baouendi, Mir, Rothschild and Zaitsev \cite{BRZ01, BMR02} provide a partial positive answer when the submanifolds are furthermore assumed to be CR. In \cite{BRZ01, BMR02}, the positive solution is obtained by approximating in the Krull topology a given formal holomorphic equivalence by a convergent one, following the spirit of Artin's approximation theorem \cite{Artin68}.  In this paper, we prove the following Artin type approximation theorem  for arbitrary smooth CR mappings of any positive codimension.

\bt\label{t:artin_gen} Let $M\subset \C^N$ be a connected real-analytic CR submanifold that is minimal at some point, $\Mp\subset \C^{\Np}$ be a real-analytic subset, and $H\colon M\to \Mp$ be a $\Cinf$-smooth CR mapping. Then there exists a dense open subset $\mc{O}\subset M$ such that for any $q\in \mc{O}$ and any positive integer $k$ there exists a germ at $q$ of a real-analytic CR mapping $H^k\colon (M,q)\to \Mp$ whose $k$-jet at $q$ agrees with that of $H$ up to order $k$.
\et

Here minimality is meant in the sense of Tumanov (see Section 2 for the precise definition). To the author's knowledge, Theorem~\ref{t:artin_gen} is the first result of its kind for mappings of positive codimension between arbitrary real-analytic submanifolds. When the target is a real-algebraic set instead of a real-analytic set, then Theorem~\ref{t:artin_gen} follows from the work of Meylan, Mir and Zaitsev \cite{MMZ03a}. Observe that Theorem~\ref{t:artin_gen} is also new even in the case $N=\Np$ since there is no rank assumption on the mapping under consideration (compare with \cite{BRZ01, BMR02, Sun08}). On the other hand, we do not know whether one may choose in Theorem~\ref{t:artin_gen} the dense open subset $\mc{O}\subset M$ to be a Zariski open subset independent of the mapping $H$. Note that when $\Mp$ is real-algebraic, such a choice is possible and follows from the main result of \cite{MMZ03a}. For more details related to Artin type approximation in CR geometry, we refer the reader to the survey paper \cite{MMZ03b}.

In this paper we shall give a rather elementary and self-contained proof of Theorem \ref{t:artin_gen}. For this, we will use several main  steps of \cite{Dam01} for which we will provide simplified proofs of the results needed for this paper. 

We will organize the paper as follows. Section \ref{s:prelim} contains some basic definitions and technical lemmas used  in Section \ref{s:properties}. In Section \ref{s:properties}, we give some elementary properties of a complex-analytic set invariantly attached to a graph of a smooth CR-mapping. The last section is devoted to the proof of Theorem \ref{t:artin_gen}.

\section{Preliminaries}\label{s:prelim}

In this section we first recall some basic definitions and  prove a lemma used in Section \ref{s:properties}. For basic background on CR analysis, we refer the reader to \cite{BERbook}. Let $M\subset\C^N$ be a real-analytic generic submanifold of codimension $d$. Let us recall that $M$ is said to be minimal at $p\in M$ if there is no germ of a real submanifold $S\subset M$ through $p$ such that the complex tangent space of $M$ at $q$ is tangent to $S$ at every $q\in S$ and $\dim_\R S < \dim_\R M$ (see \cite{BERbook}).

Following \cite{Dam01}, for a $\Cinf$-smooth CR mapping $f\colon M\to\C^{\Np}$ and for $p\in M$, we denote $\calT_p(f)$ as the germ of the smallest complex analytic set in $\C^{N+\Np}$ containing the germ of the graph of $f$ at $(p,f(p))$. The integer $\dim \calT_p(f)-N$ will be called the \textit{degree of partial analyticity} of $f$ at $p$ and denoted by $\deg_p f$. We may observe that, if $M$ is minimal at $p$, the degree of partial analyticity of $f$ at $p$ is non-negative (see Remark \ref{r:deg_pos}).

We will need the following well known reflection principle (see e.g.\ \cite{BJT85, MMZ02}).

\bp[Reflection principle]\label{p:principe_reflection}
Let $M$ be a real-analytic generic submanifold in $\C^N$ minimal at some point $p\in M$. Assume that $u\colon M\to\C$ is a continuous CR function near $p$ in $M$ and that $\T\colon\left(\C^{2N+1},\left(p,\pb,\ov{u\left(p\right)}\right)\right)\to\C$ is the germ at $\left(p,\pb,\ov{u\left(p\right)}\right)$ of a holomorphic function such that the function $F(z)=\T\left(z,\zb,\ov{u\left(z\right)}\right)$ is CR near $p$ in $M$. Then there exists a unique holomorphic extension of $F$ near $p$ in $\C^N$.
\ep

In what follows, we say that a $\Cinf$-smooth mapping $h\colon\Omega\to\C^l$, with $\Omega$ being a real manifold, is not identically zero near a point $p_0$ if the germ of $h$ at $p_0$ is non-zero, i.e. if we may find points $p$ as close as we want to $p_0$ such that $h(p)\neq 0$. 

We will also use the notion of wedge. For $p\in M$, we consider an open neighborhood $U$ of $p$ in $\C^N$ and a local defining real-analytic function $\rho\colon U\to\R^d$  of $M$ near $p$. If $\Gamma$ is an open convex cone in $\R^d$ with vertex at the origin, an open set $\calW$ of the form $\{z\in U, \rho(z,\zb)\in\Gamma\}$ is called a \emph{wedge of edge} $M$ in the direction $\Gamma$ centered at $p$.

The following result is a lemma from \cite{Dam01} for which we provide a more elementary proof.

\bl\label{l:maj_deg}
Let $M\subset\C^N$ be a real-analytic generic submanifold, minimal at $p\in M$, let $F\colon (M,p)\to\C^s$ and $u\colon (M,p)\to\C^t$ be two germs of $\Cinf$-smooth CR mappings and let $\psi\colon\left(\C^{2N+t+s}, \left(p,\pb,\ov{u\left(p\right)},F\left(p\right)\right) \right)\to\C$ be a germ at $\left(p,\pb,\ov{u\left(p\right)},F\left(p\right)\right)$ of a holomorphic function. Assume that $\psi\left(z,\zb,\ov{u(z)},F(z)\right)\equiv0$ for $z\in M$ near $p$ and that the function $(z,w)\mapsto\psi\left(z,\zb,\ov{u\left(z\right)},w\right)$ is not identically zero near $(p,F(p))$ in $M\times\C^s$. Then there exists $q\in M$ as close as we want to $p$ such that $\deg_q F<s$.
\el
\bpf
This result will be proved by induction on the integer $s$. First, we  consider the case where $s=1$.

Let $A$ be the set of points $p_1$ near $p$ in $M$ such that the holomorphic function $\C\ni w\mapsto \psi\left(p_1,\ov{p_1},\ov{u\left(p_1\right)},w\right)$ is not identically zero near $F(p_1)$ in $\C$. We may find $p_1\in A$ as close as we want to $p$. Indeed, since $\psi\left(z,\zb,\ov{u\left(z\right)},w\right)$ is not identically zero near $\left(p,F\left(p\right)\right)$ in $M\times\C$, there exists $(p_1,p^\p_1)$ as close as we want to $(p,F(p))$ such that $\psi\left(p_1,\ov{p_1},\ov{u\left(p_1\right)},p^\p_1\right)\neq0$. So, the holomorphy of $\C\ni w\mapsto\psi\left(p_1,\ov{p_1},\ov{u\left(p_1\right)},w\right)$ implies it cannot be identically zero near $F(p_1)$ in $\C$. Moreover, for $z\in A$ fixed, since the holomorphic function $\C\ni w\mapsto\psi\left(z,\ov{z},\ov{u\left(z\right)},w\right)$ doesn't vanish identically near $F(z)$ in $\C$, there exists a unique positive integer $k_z$ such that $\mderi{w}{k_z}{\psi}\left(z,\ov{z},\ov{u\left(z\right)},F\left(z\right)\right)\neq0$ and, for any integer $k<k_z$, $\mderi{w}{k}{\psi}\left(z,\ov{z},\ov{u\left(z\right)},F\left(z\right)\right)=0$. Now we fix a sufficiently small open neighborhood $V$ of $p$ in $M$ and consider the integer
$$K=\min\{k_z,z\in A\cap V\}.$$
We may pick $p_1\in A\cap V$ such that $k_{p_1}=K$, and we have
\beqar
\mderi{w}{K-1}{\psi}\left(p_1,\ov{p_1},\ov{u\left(p_1\right)},F\left(p_1\right)\right)=0,\\
\mderi{w}{K}{\psi}\left(p_1,\ov{p_1},\ov{u\left(p_1\right)},F\left(p_1\right)\right)\neq0.
\eeqar
Since $\psi$ is holomorphic near $\left(p,\pb,\ov{u\left(p\right)},F\left(p\right)\right)$ in $\C^{2N+t+1}$ by the implicit function theorem, there exists a germ at $\left(p_1,\ov{p_1},\ov{u\left(p_1\right)}\right)$ of a holomorphic function $$\T\colon\left(\C^{2N+t},\left(p_1,\ov{p_1},\ov{u\left(p_1\right)}\right)\right)\to\C$$ such that the zeros of $\mderi{w}{K-1}{\psi}$ near $\left(p_1,\ov{p_1},\ov{u\left(p_1\right)},F\left(p_1\right)\right)$ in $\C^{2N+t+1}$ are given by the equation $w=\T\left(z,\z,\nu\right)$.
On the other hand, we may observe that the function $$\mderi{w}{K-1}{\psi}\left(z,\zb,\ov{u\left(z\right)},F\left(z\right)\right)$$ is identically zero near $p_1$ in $M$. Suppose, in order to reach a contradiction, that it is false. In this case, we may find $p_2$ as close as we want to $p_1$ in $M$ such that $$\mderi{w}{K-1}{\psi}\left(p_2,\ov{p_2},\ov{u\left(p_2\right)},F\left(p_2\right)\right)\neq 0,$$ so there is a point $p_2\in A\cap V$ with $k_{p_2}< K$. This is a contradiction in view of the definition of $K$. 
So $\mderi{w}{K-1}{\psi}\left(z,\zb,\ov{u\left(z\right)},F\left(z\right)\right)\equiv0$ for $z\in M$ near $p_1$, and, from the remark on the zeros of $\mderi{w}{K-1}{\psi}$, we obtain that $F(z)=\T\left(z,\zb,\ov{u\left(z\right)}\right)$ near $p_1$ in $M$. 
But if $p_1$ is close enough to $p$, we may assume that $M$ is minimal at $p_1$ (since $M$ is real-analytic and minimal at $p$), and consequently we may apply Proposition \ref{p:principe_reflection} to obtain the existence of a holomorphic extension $\widetilde{F}$ of $F$ near $p_1$ in $\C^N$. Thus the graph of $F$ is contained, near $p_1$, in the graph of $\widetilde{F}$, which is a complex analytic set of dimension $N$. Consequently, the dimension of $\calT_{p_1}(F)$ is less than or equal to $N$. So we proved that, for an arbitrary small neighborhood $V$ of $p$ in $M$, there exists $p_1\in V$ such that $\deg_{p_1} F\leq0$. This finishes the proof of the lemma for $s=1$.

Now, we assume that the lemma holds for $s-1$, and for any $t\in\N$, any germs of CR mappings $F\colon(M,p)\to\C^{s-1}$, $u\colon (M,p)\to\C^t$ and any germ at $\left (p,\pb,\ov{u\left (p\right )},F\left (p\right )\right )$ of a holomorphic function $\psi\colon\left(\C^{2N+t+s-1}, \left(p,\pb,\ov{u\left(p\right)},F\left(p\right)\right) \right)\to\C$. Our aim is to prove the same result for $s$. We write $w=(\wp,w_s)\in\C^{s-1}\times\C$ and $F=(\Fp,F_s)\in\C^{s-1}\times\C$.

First, we consider the case where $\psi\left(z,\zb,\ov{u\left(z\right)},\Fp\left(z\right),w_s\right)\equiv0$  for $(z,w_s)\in M\times\C$ near $(p,F_s(p))$. Taking the Taylor series of $\psi$ in $w_s$ at $F_s(p)$,
\beq
\psi(z,\z,\nu,\wp,w_s)=\sum_{k\in\N}\psi_k(z,\z,\nu,\wp)\left(w_s-F_s\left(p\right)\right)^k,
\eeq
we obtain that, for any $k\in\N$, $\psi_k\left(z,\zb,\ov{u\left(z\right)},\Fp\left(z\right)\right)\equiv0$ for $z\in M$ near $p$ and that there exists $k_0$ such that $\psi_{k_0}\left(z,\zb,\ov{u\left(z\right)},\wp\right)$ doesn't vanish identically near $\left(p,\Fp\left(p\right)\right)$ in $M\times\C^{s-1}$. So, by the induction hypothesis, there exists $q\in M$ as close as we want to $p$ such that $\deg_q\Fp<s-1$, which implies $\deg_q F<s$. This completes the proof for this case.

To finish the proof, we have to consider the case where $\psi\left(z,\zb,\ov{u\left(z\right)},\Fp\left(z\right),w_s\right)$ doesn't vanish identically near $(p,F_s(p))$ in $M\times\C$. By the same method as in the case $s=1$, we show that, for point $p_1$ as close as we want to $p$ where $M$ is minimal, there exists a germ at $\left(p_1,\ov{p_1},\ov{u\left(p_1\right)},\Fp\left(p_1\right)\right)$ of a holomorphic function $$\T\colon\left(\C^{2N+t+s-1},\left(p_1,\ov{p_1},\ov{u\left(p_1\right)},\Fp\left(p_1\right)\right)\right)\to \C$$ such that 
\beq
F_s(z)=\T\left(z,\zb,\ov{u\left(z\right)},\Fp\left(z\right)\right)
\eeq
near $p_1$ in $M$.
Let $\T(z,\z,\nu,\wp)=\sum_{\a\in\N^{s-1}}\T_\a(z,\z,\nu)(\wp-\Fp(p_1))^\a$ be the Taylor series of $\T$ in $\wp$ at $\Fp(p_1)$. 

If every $\T_\a\left(z,\zb,\ov{u\left(z\right)}\right)$ is CR near $p_1$ in $M$, then by Proposition \ref{p:principe_reflection} (recall that $M$ is minimal at $p_1$) the function $M\times\C^{s-1}\ni(z,\wp)\mapsto\T\left(z,\zb,\ov{u\left(z\right)},\wp\right)\in\C$ can be holomorphically extended near $\left(p_1,\Fp\left(p_1\right)\right)$ in $\C^{N+s-1}$. We denote the extension by $\widetilde{\T}$. The graph of $F$ is contained in the graph of $\widetilde{\T}$, which is a complex submanifold of $\C^{N+s}$ of dimension $N+s-1$. This is equivalent to saying that the degree of partial analyticity of $F$ at $p_1$ is smaller than $s-1$.

If there is a multi-index $\a\in\N^{s-1}$ such that the mapping $\T_\a\left(z,\zb,\ov{u\left(z\right)}\right)$ is not CR, then there exists a vector field $\bar{L}=\sum_{j=1}^N a_j(z,\zb)\frac{\partial}{\partial \bar{z_j}}$ near $p$ in $\C^N$, where $a_1,\ldots,a_N$ are real-analytic functions near $p$, such that $\bar{L}\vert_M$ is a CR vector field and 
$$\bar{L}\left(\T\left(z,\zb,\ov{u\left(z\right)},\wp\right)\right)\nequiv0$$
near $(p_1,\Fp(p_1))$ in $M\times\C^{s-1}$. Using the chain rule, we may observe that there exists a holomorphic function $\Psi_1$ near $\left(p_1,\ov{p_1},\left.\bar{L}\left (\ov{u(z)}\right )\right \vert_{z=p_1},\Fp\left(p_1\right)\right)$ in $\C^{2N+t+s-1}$ such that 
\beq\label{chgt_ecri}
\bar{L}\left(\T\left(z,\zb,\ov{u\left(z\right)},\wp\right)\right)=\Psi_1\left(z,\zb,\bar{L}\left (\ov{u(z)}\right ),\wp\right)
\eeq
near $(p_1,\Fp(p_1))$ in $M\times\C^{s-1}$. Since $M$ is minimal at $p$, Tumanov's extension theorem (see \cite{BERbook}) implies that there exists a holomorphic extension $\tilde{u}$ of $u$ in a wedge $\calW$ of edge $M$ centered at $p$. We may assume that the mapping $\tilde{u}$  is $\Cinf$-smooth on $\calW$ up to the edge $M$. Moreover, for any $z\in \calW\cup\left (M\cap V \right )$, where $V$ is a sufficiently small neighborhood of $p$ in $\C^N$,
$$\bar{L}\left (\ov{\tilde{u}(z)}\right )=\sum_{j=1}^N a_j(z,\zb)\ov{\frac{\partial \tilde{u}}{\partial z_j}(z)}.$$
Now, for any $j\in\{1,\ldots,N$\}, $\frac{\partial\tilde{u}}{\partial z_j}$ is holomorphic in $\calW$ and $\Cinf$ up to the edge $M$. Thus, the mapping $U$ whose components are the restrictions to $M$ near $p$ of the derivatives of $\tilde{u}$ is a CR mapping near $p$ in $M$. Consequently, from the identity \eqref{chgt_ecri}, and since $\bar{L}(\ov{u(z)})=\bar{L}(\ov{\tilde{u}(z)})$, for $z\in M$ close enough to $p$, we may find a germ at $\left(p_1,\ov{p_1},\ov{U(p_1)},\Fp\left(p_1\right)\right)$ of a holomorphic function $\Psi\colon\left(\C^{2N+t+Nt+s-1},\left(p_1,\ov{p_1},\ov{U(p_1)},\Fp\left(p_1\right)\right)\right)\to\C$ such that 
$$\bar{L}\left(\T\left(z,\zb,\ov{u\left(z\right)},\wp\right)\right)=\Psi\left(z,\zb,\ov{U\left(z\right)},\Fp\left(z\right)\right)$$
near $(p_1,\Fp(p_1))$ in $M\times\C^{s-1}$. Moreover, we know that $\bar{L}\left(\T(z,\zb,\ov{u\left(z\right)},\Fp\left(z\right)\right)\equiv \bar{L}\left(F_s\left(z\right)\right)\equiv0$ near $p_1$ in $M$, since $F$ is CR; i.e. $\Psi\left(z,\zb,\ov{U\left(z\right)},\Fp\left(z\right)\right)\equiv0$ near $p_1$ in $M$. So, since we saw that $\Psi\left (z,\zb,\ov{U\left (z\right )},\wp\right )\nequiv0$ near $\left (p_1,\Fp\left (p_1\right )\right )$ in $M\times\C^{s-1}$, the induction assumption implies that the degree of partial analyticity of $\Fp$ is strictly smaller than $s-1$ for points in $M$ as close as we want to $p_1$ and therefore as close to $p$. This finishes the proof of Lemma \ref{l:maj_deg}.
\epf

As in \cite{Dam01}, one gets from Lemma \ref{l:maj_deg} the following result.

\bl\label{l:maj_deg_bis}
Let $M\subset\C^N$ be a real-analytic generic submanifold, minimal at $p\in M$, let $F\colon (M,p) \to \C^s$ be a germ of a CR mapping and let $$\psi\colon\left(\C^{2N+2s}, \left(p,\pb,\ov{F\left(p\right)},F\left(p\right)\right) \right)\to\C$$ be a germ  at $\left(p,\pb,\ov{F\left(p\right)},F\left(p\right)\right)$ of a holomorphic function. Assume that $\psi\left(z,\zb,\ov{F(z)},F(z)\right)\equiv0$ for $z\in M$ near $p$ and that the function $(z,v,w)\mapsto\psi(z,\zb,v,w)$ is not identically zero near $\left(p,\ov{F\left(p\right)},F\left(p\right)\right)$ in $M\times\C^{2s}$. Then, there exists $q\in M$ as close as we want to $p$ such that $\deg_q F<s$.
\el

\bpf
First, we assume that $\psi\left(z,\zb,\ov{F\left(z\right)},w\right)$ does not vanish identically near $(p,F(p))$ in $M\times\C^s$. Since $\psi\left(z,\zb,\ov{F\left(z\right)},F(z)\right)\equiv0$ for $z\in M$ near $p$, we may apply Lemma \ref{l:maj_deg} and deduce that the degree of partial analyticity of $F$ at $q$ is  strictly smaller than $s$ for $q$ arbitrarily close to $p$.

Now we treat the case where $\psi\left(z,\zb,\ov{F\left(z\right)},w\right)\equiv0$ for $(z,w)\in M\times\C^s$ near $(p,F(p))$. For this we consider the Taylor series of $\psi$ in $w$ at $F(p)$:
\beq
\psi(z,\z,v,w)=\sum_{\a\in\N^s}\psi_\a(z,\z,v)\left(w-F\left(p\right)\right)^\a.
\eeq
The assumption implies that, for any $\a\in\N^s$, $\psi_\a\left(z,\zb,\ov{F\left(z\right)}\right)\equiv 0$ near $p$ in $M$. However, by assumption, there is a multi-index $\a_0$ such that $\psi_{\a_0}\left(z,\zb,v\right)\nequiv 0$ near $\left (p,\ov{F\left (p\right )}\right )$ in $M\times\C^s$, and Lemma \ref{l:maj_deg} gives the desired result.
\epf

\section{Properties of $\calT_p(H)$}\label{s:properties}
In this section, we fix a real-analytic generic submanifold $M\subset\C^N$ and a $\Cinf$-smooth CR mapping $H\colon M\to\Mp$ on $M$ with values in a real-analytic set $\Mp\subset\C^{\Np}$. We shall give some properties of the degree of partial analyticity of $H$ and of the complex analytic set $\calT_p(H)$ for $p\in M$. All the results of this section can be found in \cite{Dam01}, but since the proofs we shall give are rather elementary compared to \cite{Dam01}, we include them in this note for completeness. 

The following lemma is a direct consequence of the boundary uniqueness theorem.

\bl\label{l:prolon&graph}
Let $M$, $\Mp$ and $H$ be as above. If $H$ admits a holomorphic extension $\Htilde$ on a wedge of edge $M$ centered at $p$, then the graph of $\Htilde$ near $\left(p,H(p)\right)$ is contained in $\calT_p(H)$.
\el

\brm\label{r:deg_pos}
If $M$ is minimal at $p\in M$, Tumanov's extension theorem implies that, in the above setting, there is a unique extension $\Htilde$ of $H$ holomorphic on a wedge of edge $M$ centered at $p$. Thus $\calT_p(H)$ contains the graph of $\Htilde$ near $(p,H(p))$. This implies that the dimension of $\calT_p(H)$ at $p$ is greater than $N$, i.e. that the degree of partial analyticity of $H$ at $p$ is non-negative.
\erm

The two following lemmas describe the regular points of the complex analytic set $\calT_p(H)$.

\bl\label{l:pt_reg}
In the above setting, assume that the degree of partial analyticity of $H$ is constant near some point $p$ in $M$. Then there exists an open neighborhood $U^1_p$ of $p$ in $M$ such that the set $\Sigma^1_p\subset M\cap U^1_p$ of points $q\in U^1_p$ for which $\calT_p(H)$ is not regular at $\left(q,H(q)\right)$ is a closed set with empty interior.
\el

\bpf
We may find an open neighborhood $U^1_p$ of $p$ in $M$ on which $\calT_p(H)$ is a complex analytic set and on which the degree of partial analyticity of $H$ is constant equal to $s$. Let $\Sigma^1_p\subset M\cap U^1_p$ be the set of points $q\in U^1_p$ for which $\calT_p(H)$ is not regular at $\left(q,H(q)\right)$. In view of the classical definition of regular points of a complex analytic set, $\Sigma^1_p$ is a closed subset of $M\cap U^1_p$. To prove that its interior is empty, assume by contradiction that we may find an open subset $V$ of $M$ contained in $\Sigma^1_p$. Thus, for any $q\in V$, the graph of $H$ near $(q,H(q))$ is contained in the set of the singular points of $\calT_p(H)$, which is a complex analytic set of dimension strictly smaller than $N+s$. So, the dimension of $\calT_q(H)$ is also strictly smaller than $N+s$ for any $q\in V$. This is impossible, since the degree of partial analyticity of $H$ is constant equal to $s$ on $V$.
\epf

\bl\label{l:equ_def}
In the above setting, assume that $M$ is minimal at $p\in M$, that the degree of partial analyticity of $H$ is constant equal to $s$ near $p$ in $M$, and write $t=\Np-s$. Then there are an open neighborhood $U^2_p\subset U^1_p$ of $p$ in $M$ (with $U^1_p$ given by Lemma \emph{\ref{l:pt_reg}}) and a closed set with empty interior $\Sigma^2_p\subset M\cap U^2_p$ such that, for any $q\in U^2_p\setminus\Sigma^2_p$, there are holomorphic coordinates $(u^\p,v^\p)\in\C^s\times\C^t$ near $H(q)$ for which $H=(F,G)\in\C^s\times\C^t$, and a germ at $(q,F(q))$ of a holomorphic mapping $T_q\colon\left(\C^{N+s},(q,F(q))\right)\to\C^t$ such that $\calT_p(H)$ is given near $(q,H(q))$ by the equation $v^\p=T_q\left(z,u^\p\right)$. In particular $\calT_p(H)$ is regular at $(q,H(q))$.
\el

\bpf
Since $M$ is minimal at $p$ and $H$ is CR on $M$, Tumanov's extension theorem implies that there exists a holomorphic extension $\Htilde$ of $H$ in a wedge $\calW$ of edge $M$ centered at $p$ which is $\Cinf$-smooth up to the edge. By Lemma \ref{l:prolon&graph} the graph of $\Htilde$ is contained in $\calT_p(H)$ near $(p,H(p))$. It means that we may choose an open neighborhood $\Delta$ of $p$ in $\C^N$ such that $(z,\Htilde(z))\in \calT_p(H)$ for any $z\in\Delta\cap\calW$.

On the other hand, if $U^1_p$ is the open neighborhood of $p$ given by Lemma \ref{l:pt_reg}, we define $U^2_p=U^1_p\cap\Delta$ and $\Sigma^2_p=\Sigma^1_p\cap U^2_p$. Now, for a fixed point $q\in U^2_p\setminus \Sigma^2_p$, $\calT_p(H)$ is regular at $(q,H(q))$. Moreover, since the degree of partial analyticity of $H$ is constant equal to $s$ on $U^1_p$, the dimension of $\calT_p(H)$ at $(q,H(q))$ is $N+s$. So, there exist an open neighborhood $\mc{U}\subset U^2_p$ of $q$ in $\C^{N}$, an open neighborhood $\mc{V}$ of $H(q)$ in $\C^{\Np}$ and a holomorphic mapping $f\colon\mc{U}\times\mc{V}\to\C^t$ of rank $t$ at $(q,H(q))$ such that 
\beq
\calT_p(H)=\{(z,\zp)\in \mc{U}\times\mc{V},\ f\left (z,\zp\right )=0 \}.
\eeq
Thus, for $z\in\mc{U}^\p=\mc{U}\cap\calW\cap\Htilde^{-1}(\mc{V})$, we have 
\beq
f\left (z,\Htilde\left (z\right )\right )=0.
\eeq
So $\deri{z}{f}\left (z,\Htilde\left (z\right )\right )+\deri{\zp}{f}\left (z,\Htilde\left (z\right )\right )\cdot\deri{z}{\Htilde}\left (z\right )\equiv 0$ in $\mc{U}^\p$. Since the mapping $\Htilde$ is $\Cinf$-smooth up to the edge $M$, we have the following identity: 
\beq\label{equ_rk}
\deri{z}{f}\left (z,\Htilde\left (z\right )\right )+\deri{\zp}{f}\left (z,\Htilde\left (z\right )\right )\cdot\deri{z}{\Htilde}\left (z\right )\equiv 0
\eeq
on $M\cap\mc{U}^\p$. But $f$ is of rank $t$ at $(q,H(q))$ and the identity $\eqref{equ_rk}$ shows that the columns of $\deri{z}{f}\left (q,H\left (q\right )\right )$ are a linear combination of the columns of $\deri{\zp}{f}\left (q,H\left (q\right )\right )$. Consequently the rank of  $\deri{\zp}{f}(q,H(q))$ is $t$. So, by the implicit function theorem there exist holomorphic coordinates $(u^\p,v^\p)\in\C^s\times \C^t$ for which $H=(F,G)\in\C^s\times\C^t$ and a holomorphic mapping $T_q$ near $(q,F(q))$ in $\C^{N+s}$ such that the zeros of $f$ are given by points of the form $\left (z,u^\p,T_q\left (z,u^\p\right )\right )\in\C^N\times\C^s\times\C^t$.
\epf

For the next result, we will denote $\pi\colon\C^N\times\C^{\Np}\to\C^N$ and $\pi^\p\colon\C^N\times\C^{\Np}\to\C^{\Np}$ as the canonical projections.
 
\bl\label{l:projectionbis}
In the above setting, assume that $M$ is minimal at $p\in M$ and that the degree of partial analyticity of $H$ is constant equal to $s$ near $p$ in $M$. Then there are an open neighborhood $U^3_p\subset U^2_p$ of $p$ in $M$ (with $U^2_p$ given by Lemma \emph{\ref{l:equ_def}}) and a closed set with empty interior $\Sigma^3_p\subset M\cap U^3_p$ such that, for any $q\in U^3_p\setminus\Sigma^3_p$, there is a neighborhood $\Omega_q$ of $(q,H(q))$ in $\C^{N+\Np}$ satisfying
$$\pi^\p\left(\calT_p\left(H\right)\vert_{M\times\C^{\Np}}\cap\Omega_q\right)\subset\Mp.$$
\el

\bpf
Let $\rp\colon\C^{\Np}\to\R^{d^\p}$ be a local defining real-analytic function of $\Mp$ near $H(p)$. Since $H(M)\subset\Mp$, we have the mapping identity 
\beq\label{inclus}
\rp\left (H(z),\ov{H(z)}\right )\equiv0
\eeq
for $z\in M$ near $p$.

We consider $U^2_p$ and $\Sigma^2_p$ respectively as the open neighborhood of $p$ in $M$ and the closed set with empty interior given by Lemma \ref{l:equ_def}. Let $U^3_p\subset U^2_p$ be a sufficiently small connected neighborhood of $p$ in $M$.
Since $M$ is real-analytic and minimal at $p$, the set $\Sigma_1$ of points in $U^3_p$ where $M$ is not minimal is a closed set with empty interior. Thus, $\Sigma^3_p=\left (\Sigma^2_p\cap U^3_p\right )\cup \Sigma_1$ is a closed set with empty interior. Moreover, for any $q\in U^3_p\setminus\Sigma^3_p$, $M$ is minimal at $q$, and we may find holomorphic coordinates $(u^\p,v^\p)\in\C^s\times\C^t$ near $H(q)$ for which $H=(F,G)\in\C^s\times\C^t$, and a germ at $(q,F(q))$ of a holomorphic mapping $T_q\colon\left(\C^{N+s},(q,F(q))\right)\to\C^t$ such that $\calT_p(H)$ is given near $(q,H(q))$ by the equation $v^\p=T_q\left(z,u^\p\right)$. If $U^3_p$ is chosen small enough we may assume that the graph of $H$ is contained near $q\in U^3_p\setminus\Sigma^3_p$ in $\calT_p(H)$; i.e. we have the mapping identity 
\beq\label{equ_H_1}
G(z)=T_q(z,F(z))
\eeq
for $z\in M$ near $q$. We may also assume that the mapping defined by 
\beq\label{def_psi}
\Psi(z,\z,\nu^\p,u^\p)=\rp\left(u^\p,T_q\left(z,u^\p\right),\nu^\p,\ov{T_q\left(\ov{\z},\ov{\nu^\p}\right)}\right)
\eeq
is holomorphic near $(q,\bar{q},\ov{F(q)},F(q))$ in $\C^{2N+2s}$ for any $q\in U^3_p\setminus\Sigma^3_p$.

From $\eqref{equ_H_1}$ and $\eqref{inclus}$, we obtain
\beq\label{inclus_2}
\rp\left (F\left (z\right ),T_q\left (z,F\left (z\right )\right ),\ov{F(z)},\ov{T_q\left (z,F(z)\right )}\right )\equiv0
\eeq
for $z\in M$ near $q$.
Thus, by $\eqref{inclus_2}$, we have
\beq\label{equ1_psi}
\Psi\left (z,\zb,\ov{F(z)},F(z)\right )\equiv0
\eeq
for $z\in M$ near $q$. Since the degree of analyticity of $H$ is constant equal to $s$ on $U^3_p$, from $\eqref{equ_H_1}$, we deduce that the degree of partial analyticity of $F$ is constant equal to $s$ on $U^3_p$. So, the identity $\eqref{equ1_psi}$ and Lemma \ref{l:maj_deg_bis} (recall that $M$ is minimal at $q$) imply that 
\beq
\Psi(z,\zb,\nu^\p,u^\p)\equiv 0
\eeq
for $(z,\nu^\p,u^\p)\in M\times\C^{2s}$ near $(q,\ov{F(q)},F(q))$. This identity is equivalent to 
\beq\label{eq_final}
\rp\left (u^\p,T_q(z,u^\p),\nu^\p,\ov{T_q(z,\ov{\nu^\p})}\right )\equiv 0
\eeq
for $(z,\nu^\p,u^\p)\in M\times\C^{2s}$ near $\left (q,\ov{F(q)},F(q)\right )$. 
So taking $\nu^\p=\ov{u^\p}$ in \eqref{eq_final}, we obtain that $(z,u^\p,T_q(z,u^\p))\in\Mp$ as soon as $z$ is close enough to $q$ in $M$ and $u^\p$ is close enough to $F(q)$ in $\C^s$. This finishes the proof of the lemma.
\epf

We conclude this section with a lemma which makes use of the upper semi-continuity of the partial analyticity degree; we leave the details to the reader.

\bl\label{l:degr_cst}
In the above setting, there exists a closed set with empty interior $\Sigma_2\subset M$ such that the degree of partial analyticity of $H$ is constant on each connected component of $M\setminus \Sigma_2$.
\el

\section{Proof of Theorem \ref{t:artin_gen}}\label{s:proofs}

In this section, we keep the assumptions and the notation of Section \ref{s:properties} i.e we consider a real-analytic generic submanifold $M$ in $\C^N$ and a $\Cinf$-smooth CR mapping $H\colon M\to\Mp$ where $\Mp\subset\C^{\Np}$ is a real-analytic set. Theorem \ref{t:artin_gen} will be a consequence of the following result.

\bp\label{t:artin_vois_pt}
In the above setting, assume that $M$ is minimal at $p\in M$ and that the degree of partial analyticity of $H$ is constant equal to $s$ near $p$ in $M$. Let $U^3_p$ and $\Sigma^3_p\subset M\cap U^3_p$ be, respectively, the open neighborhood of $p$ in $M$ and the closed set with empty interior given by Lemma \ref{l:projectionbis}. Then for any $q\in U^3_p\setminus\Sigma^3_p$ and for any positive integer $k$ there exists a germ at $q$ of a real-analytic CR mapping $H^k\colon\left(M,q\right)\to\Mp$ whose $k$-jet at $q$ agrees with that of $H$ up to order k.
\ep
\bpf
For any $q\in U^3_p\setminus\Sigma^3_p$, we may choose holomorphic coordinates  $\zp=(u^\p,v^\p)\in\C^s\times\C^t$ and a germ at $(q,F(q))$ of a holomorphic mapping $T_q\colon\left(\C^{N+s},(q,F(q))\right)\to\C^t$ such that $\calT_p(H)$ is given by the equation $v^\p=T_q(z,u^\p)$ near $(q,H(q))$. We may also find a neighborhood $\Omega_q$ of $(q,H(q))$ in $\C^{N+\Np}$ such that $\pi^\p\left(\calT_p\left(H\right)\vert_M\cap\Omega_q\right)\subset\Mp.$

We fix a point $q\in U^3_p\setminus\Sigma^3_p$ and a positive integer $k$. Since the graph of $H$ is contained in $\calT_p(H)$ near $\left(q,H\left(q\right)\right)$, we have the following mapping identity:
\beq\label{equ_H}
G(z)=T_q(z,F(z))
\eeq 
near $q$ in $M$.
Let $F^k$ be the $k$-th order Taylor polynomial of $F$ at $q$ (that is holomorphic since $H$ is CR; see \cite{BERbook}). We define the holomorphic mapping  $G^k$ near $q$ in $\C^N$ by setting 
\beq\label{def_Gk}
G^k(z)=T_q\left(z,F^k\left(z\right)\right).
\eeq
Thus $H^k=\left (F^k,G^k\right )$ is a holomorphic mapping near $q$ in $\C^N$ with values in $\C^{\Np}$. Moreover, by definition of $H^k$ and from $\eqref{equ_H}$, the $k$ first derivatives of $H^k$ at $q$ coincide with that of $H$. So to complete the proof of Proposition \ref{t:artin_vois_pt}, we have to show that $H^k$ sends $M$ into $\Mp$ near $q$. From the definition of $H^k$ and the local defining equation of $\calT_p(H)$ near $\left(q,H\left(q\right)\right)$ in $\C^{N+\Np}$, we obtain that there exists an open neighborhood $\widetilde{\Omega_q}$ of $\left(q,H\left(q\right)\right)$ in $\C^{N+\Np}$ such that 

$$\mc{G}_{H^k}\cap \left(M\times\C^{\Np}\right)\cap\widetilde{\Omega_q}\subset\calT_p(H)\vert_{M\times\C^{\Np}}\cap\Omega_q,$$
where $\mc{G}_{H^k}$ is the graph of the mapping $H^k$. Since $\pi^\p\left(\calT_p\left(H\right)\vert_M\cap\Omega_q\right)\subset\Mp$ by Lemma~\ref{l:projectionbis}, the conclusion of the proposition follows.
\epf

Now, we are able to prove our main result, Theorem \ref{t:artin_gen}.

\bpf[Proof of Theorem \ref{t:artin_gen}] We first note that Theorem~\ref{t:artin_gen} holds in the case where $N=1$ or $N'=1$. We may therefore assume that $N,N'\geq 2$.

We first treat the case where $M$ is generic. Since $M$ is real-analytic and connected, there exists a real-analytic subvariety $\Sigma_1$ of $M$ such that $M$ is minimal at each point $p\in M\setminus \Sigma_1$. Moreover, from Lemma \ref{l:degr_cst}, there exists a closed set with empty interior $\Sigma_2\subset M$ such that the degree of partial analyticity of $H$ is constant on each connected components of $M\setminus\Sigma_2$. Thus, $\Sigma=\Sigma_1\cup\Sigma_2$ is a closed subset of $M$ with empty interior. Fix a point $p\in M\setminus\Sigma$; by Proposition \ref{t:artin_vois_pt}, we may find an open subset $U^3_p$ of $p$ in $M$ and a closed set with empty interior $\Sigma^3_p\subset M\cap U^3_p$ such that, for every $q\in U^3_p\setminus\Sigma^3_p$, the conclusion of Theorem \ref{t:artin_gen} holds at $q$. Thus, the set $\mc{O}=\bigcup_{p\in M\setminus \Sigma} \left(U^3_p\setminus \Sigma^3_p\right)$ does the job, and this finishes the proof of Theorem \ref{t:artin_gen} for the generic case.

If $M$ is not generic, for any $p\in M\setminus\Sigma_1$ (where $\Sigma_1$ again denotes the set of nonminimal points of $M$) we may assume, thanks to a local holomorphic change near $p$, that $M=\widetilde{M}_p\times\{0\}\subset \C^{N-r_1}_{z_1}\times\C^{r_1}_{z_2}$, where $r_1$ is a non-negative integer and $\widetilde{M}_p$ is a connected real-analytic generic submanifold which is minimal (see \cite{BERbook}). From the generic case treated above, there exists a dense open subset  $\widetilde{\mc{O}}_p\subset\widetilde{M}_p$ such that, for any non-negative integer $k$ and any $\tilde{q}\in\widetilde{\mc{O}}_p$, there exists a germ at $\tilde{q}$ of a real-analytic CR mapping $H_1^k\colon (\widetilde{M}_p,q)\to M'$ whose $k$-jet at $\tilde{q}$ agrees with that of $\widetilde{M}_p \ni z_1\mapsto H(z_1,0)$. Since $\bigcup_{p\in M\setminus\Sigma_1}\left (\widetilde{\mc{O}}_p\times\{0\}\right )$ is a dense open subset of $M$, the proof of Theorem \ref{t:artin_gen} is complete.
\epf

\end{document}